\documentclass[12pt, reqno]{amsart}
\usepackage{amsmath, amsthm, amscd, amsfonts, amssymb, graphicx, color}
\usepackage[bookmarksnumbered, colorlinks, plainpages]{hyperref}
\hypersetup{colorlinks=true,linkcolor=red, anchorcolor=green, citecolor=cyan, urlcolor=red, filecolor=magenta, pdftoolbar=true}

\begin{document}
\title{What to do to have your paper rejected!}
\author{Mohammad Sal Moslehian$^1$ and Rahim Zaare-Nahandi$^2$}

\address{$^1$Department of Pure Mathematics, Ferdowsi University of
Mashhad, P.O. Box 1159, Mashhad 91775, Iran;\newline
Department of Mathematics and Computer Science,
Karlstad University, SE-65188 Karlstad, Sweden.}
\email{moslehian@um.ac.ir, mohammad.moslehian@kau.se}
\address{$^2$University of Tehran, School of Mathematics, Statistics and Computer Sciences, Iran}
\email{rahimzn@ut.ac.ir}



\begin{abstract}
We aim to highlight certain points and considerations for graduate students and young researchers, which should be avoided in submissions to good research journals. Observing these remarks could substantially decrease the probability of rejections of papers.
 \end{abstract}
\maketitle

Do you really want to have your paper rejected? Here is some practical advice, based on the authors' experience as reviewers and editors of international journals. Following this ``advice'' will considerably increase the probability that your paper will not be accepted by any respectable journal.

\section{mathematical content}

\noindent $\bullet~$\textbf{Write a two lines abstract expressing some very general facts.} For example,
{\footnotesize \begin{center}
``In this paper we investigate some general inequalities and
present some interesting applications.''
\end{center}
or
\begin{center}
``This paper is devoted to the study of some second order differential equations and their solutions. Our results extend some known results in the literature.''
\end{center}}\vspace{0.5cm}

\noindent $\bullet~$\textbf{Do not define the basic concepts considered in your paper.}

Then you may receive a message with the following context:
\begin{quotation}
\textit{``The exposition of the paper is poor. For example, the author does
not even provide definitions of the basic objects considered in the
paper.''}
\end{quotation}\vspace{0.5cm}

\noindent $\bullet~$\textbf{Do not provide motivations for the subject of your paper, or
ignore describing your methods.}

Then it is natural to expect the following comment from the editor:
\begin{quotation}
\textit{``While your work appears to be mathematically correct, it is not
clear what impact such results have. Publications in this journal
require clear reasons for the interest in the subject and
development of new techniques. Accordingly, your paper should be
rejected.''}
\end{quotation}\vspace{0.5cm}

\noindent $\bullet~$\textbf{Prove your results based on many results of papers which are not yet accepted
for publication or are not yet confirmed.}

Then you will receive a message like the following:
\begin{quotation}
\textit{``Some of your results are based on unpublished results which are not available on internet and could
not be confirmed by our referees. Accordingly, we suggest the author to withdraw the current paper and resubmit it after getting acceptance of the previously submitted paper(s) on which your results are based.''}
\end{quotation}\vspace{0.5cm}

\noindent $\bullet~$\textbf{Write your article just based on trivial generalizations.} For
example,
{\footnotesize \begin{center} add or reduce a parameter in an equation, or, if a property is proved for two elements, state and prove it for three elements (and plan to do it for $n$ elements in the next paper!).\end{center}}

Then you should expect a reply with the following context:
\begin{quotation}
\textit{``Trivial operations such as changing or adding a parameter on
someone else's paper do not generally lead to an original work. Many
of the resulting statements are straightforward. The readership for such a paper is
usually very limited.''}
\end{quotation}
or\\
\begin{quotation}
\textit{``The authors simply extend a known inequality involving some double
integrals to another inequality for triple integrals. Neither
serious nor new techniques have been presented. Probably the authors will next try
to publish a paper for multiple integrals!''}
\end{quotation}
or
\begin{quotation}
\textit{``Most parts of this article are well-known, and the notion of
$k$-ring does not seem to give better proofs compared to the
standard ones. The article seems to be far below the standards of the journal.''}
\end{quotation}\vspace{0.5cm}

\noindent $\bullet~$\textbf{Take a famous open conjecture, suppose one variant of answering it as a hypothesis of your theorems and build up a stack of derived facts''.}\vspace{0.5cm}

\noindent $\bullet~$\textbf{Concentrate on a worthless problem and make it your research
topic.}

Then more likely you will face with the following reply from
the editor:
\begin{quotation}
\textit{``The subject of this paper is away from the main stream of research in
mathematics. There are a few readers who are interested in such a
topic. I suggest the rejection of the paper.''}
\end{quotation}\vspace{0.5cm}

\noindent $\bullet~$\textbf{Write a paper with short and trivial proofs.} For example,
{\footnotesize \begin{center} write a 10 pages paper in which 6 pages are devoted to the introduction, 2 pages are for references, and only 2 pages are considered for the main results while more than half of these 2 pages are filled with some lemma's from other papers.
\end{center}}\vspace{0.5cm}

\noindent $\bullet~$\textbf{Omit serious parts of proofs, where there are, e.g., some gaps.} For example, write
{\footnotesize The proof is trivial and is left to the reader.}\vspace{0.5cm}

\noindent $\bullet~$\textbf{Write a paper where the statement of the theorems are too longer than their proofs.}\vspace{0.5cm}

\noindent $\bullet~$\textbf{Do not provide any nontrivial objective or abstract examples for
the concepts you define.}  For example,
{\footnotesize \begin{center}introduce ``probabilistic
non-Archimedean Jordan $CQ^*$-algebra'', which has no interesting impact beyond a definition. Or, investigate
on a property while the set of mathematical objects satisfying that
property is empty.\end{center}}

You may then receive a report as follows.
\begin{quotation}
\textit{``The author creates new terms which are close to terms in the literature (and indeed redefines existing terminology) and gives no motivation for the new or changed definition and no examples which differentiate between previous terms and new ones. The paper should be rejected.}
\end{quotation}\vspace{0.5cm}

\noindent $\bullet~$\textbf{Avoid any illustrating explanation of your results, offer only highly abstract thoughts.} \vspace{0.5cm}

\noindent $\bullet~$\textbf{Whenever you face a difficulty on proofs, or you needed a
property, immediately add it to your hypothesis. Fill your result with
so many hypotheses that no nontrivial mathematical object could
satisfy your hypotheses.}

Most likely what you will receive will be:
\begin{quotation}
\textit{''The assumptions are too strong and not interesting. The whole
paper seems to be artificial. I suggest the rejection of the paper.``}
\end{quotation}\vspace{0.5cm}

\noindent $\bullet~$\textbf{End your paper with a lemma without any application.}

\section{Style oriented}

\noindent $\bullet~$\textbf{Choose a lengthy, general or ambiguous title for your paper} such as
{\footnotesize  \begin{center}
``A Result in Group Theory''
\end{center}
or
\begin{center}
``On Banach Algebras''.
\end{center}}\vspace{0.5cm}

\noindent $\bullet~$\textbf{In the ``Style File'' ignore to delete ``Insert your abstract here.''
and write your abstract in the sequel of this statement.} For example, the
resulting text is:
{\footnotesize \begin{center}``Insert your abstract here. In this paper we determine the fool's solitaire number for the join of graphs ...''\end{center}}\vspace{0.5cm}

\noindent $\bullet~$\textbf{Present your article in a single section without giving
introduction and literature review of the subject.}\vspace{0.5cm}

\noindent $\bullet~$\textbf{Begin each section with a definition or theorem without any brief
statement for every section.} For example:\\
{\footnotesize \begin{center}
``2. Main Results\\
~Theorem 2.1. Let R be a commutative ring. ...''
\end{center}}\vspace{0.5cm}

\noindent $\bullet~$\textbf{Write your results one after another without any explanation. For
example, do not explain what Theorem A is about, or how it is
related to other results of your paper. Even write your definitions
consecutively with no interruptions.} For example:\\
{\footnotesize \begin{center}
``Definition 1. ...\\
~Definition 2. ...\\
~Definition 3. ...\\
~Definition 4. ...''\\
\end{center}}\vspace{0.5cm}

\noindent $\bullet~$\textbf{Don't give a thought to the length or writing style expected by the journal to which you send it.}\vspace{0.5cm}

\noindent $\bullet~$\textbf{Do not polish your paper and do not give your
paper to a friend for correcting typos and grammatical errors, but instead, immediately submit it to a journal.} For example write:
{\footnotesize \begin{center}
''Let A is a comutative ring and a is belong to A that is nilpotent element.``
\end{center}}
You will then receive an answer as follows
\begin{quotation}
\textit{``The presentation of the paper is unacceptable, several misprints
and typos can be detected. The paper has poor English. It is impossible for me to
scientifically follow discussions and write a review
on this paper.''}
\end{quotation}
or\\
\begin{quotation}
\textit{``It seems that the author first wrote the paper in his native language and after that translated it word by word to English with the help of a translator. It is not clear what is the mathematical meaning.''}
\end{quotation}\vspace{0.5cm}

\noindent $\bullet~$\textbf{Write a weak paper and submit it to a journal of high quality.}

You will then receive a message like:
\begin{quotation}
\textit{``The results of this paper are not substantial to merit
its publication in ...''}
\end{quotation}
or\\
\begin{quotation}
\textit{``This paper does not fulfill the general quality and novelty which
normally characterize papers published in this journal.''}
\end{quotation}\vspace{0.5cm}

\noindent $\bullet~$\textbf{Submit a paper in a subject such as algebra to a journal devoted
to geometry.}

Then you will typically receive the following response:
\begin{quotation}
\textit{``I regret to inform you that your paper is not in the scope of this
journal. You may send it to a journal matched with the topic of your
article.''}
\end{quotation}\vspace{0.5cm}

\noindent $\bullet~$\textbf{Submit your paper to a journal whose editors' interests are far from the subject of your paper.}\vspace{0.5cm}

\noindent $\bullet~$\textbf{Give incorrect Mathematical Subject Classifications and suggest non-expert referees as potential reviewers for your paper.}\vspace{0.5cm}

\noindent $\bullet~$\textbf{Provide a paper with so many irrelevant and unnecessary references.} For example,
{\footnotesize \begin{center} write a 6 pages paper with 36 references.\end{center}}

Then, you may receive a message from the journal authorities similar to the
following one:
\begin{quotation}
\textit{``Your short paper has too many references. The number of
papers/books in the list of references is expected to be about 1.5
times the number of pages of the paper and the number of
self-citations could be at most a quarter of all references unless
they have been used out of literature review and preliminaries.''}
\end{quotation}\vspace{0.5cm}

\noindent $\bullet~$\textbf{Use different formats for your references, or give incomplete
references, or provide inaccurate information.} For example:\\
{\footnotesize \begin{center}
``[1] W.B. Arveson, \textit{$C\sp *$-algebras and
numerical linear algebra}, J. Funct. Anal. \textbf{122} (1994), no.
2, 333--360.\newline
[2] Bhatia, R., Matrix Theory, (Graduate text in Mathematics) Springer Verlag, New York, 1997. \newline
[3] A. B\"ottcher, A.V. Chithra and M.N.N. Namboodiri:
Approximation of Approximation Numbers by Truncation. \textit{J.
Integr. Equ. Oper. Theory} 2001, 39, 387-395.''
\end{center}}

You may then receive the following comment:
\begin{quotation}
\textit{``The author does not take care of references. He/She should unify
the references. The paper cannot be considered in the present
form.''}
\end{quotation}


\section{Ethics}

\noindent $\bullet~$\textbf{In your abstract, mention that you have generalized the results
of a low level paper, or unreasonably acknowledge a well-known mathematician  at the end of your, e.g. 4 pages, paper.}\vspace{0.5cm}

\noindent $\bullet~$\textbf{``Insert sections and materials from another article in your paper, while to avoid accusation for plagiarism, make some small changes.''} For instance,\\
{\footnotesize if your paper is computational, use a different but somehow similar boundary value. Or, simply change the symbols, e.g. use $\alpha$-derivation instead of $\sigma$-derivation.}

Then you could expect receiving a letter with the following context:
\begin{quotation}
\textit{``All results are either known or an adaptation of results appearing
in the references. This translation to a ``new scope'' does not bring
enough interest to be published, since no new problems are solved with this approach.''}
\end{quotation}\vspace{0.5cm}

\noindent $\bullet~$\textbf{Implement small modifications on another paper to provide a new
paper.} For example,
{\footnotesize \begin{center} if the original paper is on a linear operator
$T$, replace $T$ with its adjoint $T^*$ throughout the paper and make straightforward modifications on the proofs of the original paper for $T^*$ instead of $T$.
\end{center}}\vspace{0.5cm}

\noindent $\bullet~$\textbf{Introduce your own notation for things that have standard notation.}

Then the referee may reply as follows.
\begin{quotation}
\textit{``I spent hours interpreting what is written only to discover that the actual argument or result is standard. The paper is not suitable for publication.}
\end{quotation}\vspace{0.5cm}

\noindent $\bullet~$\textbf{Call your own results interesting or well-known. Or, define an
object or a construction and put your name on this object or
construction.}

You will then receive a reply with the following context:
\begin{quotation}
\textit{``The terminology of the paper is unusual. It is not expected from
an author to put his/her name on a construction or a theorem.''}
\end{quotation}\vspace{0.5cm}

\noindent $\bullet~$\textbf{Submit your paper to more than one journal at the same time.}

You may receive the following from the editorial board
\begin{quotation}
\textit{``Your action on submitting your manuscript to two journals at the same time is unethical. According to the journal policy we have to put your name in the black list. This means that we no longer consider any submission of you for possible publication in our journal for 5 years.}
\end{quotation}

\noindent $\bullet~$\textbf{In your introduction, when referring to the latest results on the
subject, only mention your own work and refer to them.} For example, write
{\footnotesize \begin{center}``The author investigates ..., see author's papers [1, 2, 3, 4, 5, 6, 7]''.\end{center}}
Then, you may receive a message from the journal's editor similar to:
\begin{quotation}
\textit{``The literature review of the paper is too poor. The author needs
to search MathSciNet, Zentralblatt Math and Google Scholar, to find other new contributions in
the subject of the paper and cite the most important items
properly.''}
\end{quotation}\vspace{0.5cm}

\noindent $\bullet~$\textbf{Send consecutive messages to the editor of the journal for unusual
 expectations such as speeding up the refereeing procedure of your article. And, that receiving an acceptance for your paper in a short time is extremely vital for you.}

Maybe the most polite response you would receive is
\begin{quotation}
\textit{``Your paper is still under review. The handling editor of your
paper has tried to provide a report as early as possible. A peer
review essentially depends on the referee and sometimes is out of the control of the
editors. So please have more patience. However, if you feel it is not possible for you to wait, you may withdraw
your paper and after our confirmation, submit it somewhere else.''}
\end{quotation}
\bigskip
\textbf{Acknowledgement.} The authors would like to thank most editors of Banach J. Math. Anal. and Ann. Funct. Anal. for their useful comments improving the paper.

\end{document}